\documentstyle[11pt,fleqn]{article}
                      
\title{{\Large SUBMODULES OF SPECHT MODULES FOR WEYL GROUPS}} 

\author{$by$ SA\.{I}T ~HALICIO{\u G}LU} 
\date{ }
\renewcommand{\baselinestretch}{1.20} 
\begin{document}
\maketitle
\renewcommand{\baselinestretch}{0.90}

\noindent
{\footnotesize The construction of all
 irreducible modules of the symmetric groups over an arbitrary field 
 which reduce to Specht modules in the case of fields of characteristic
  zero is given by G.D.James. Hal\i c\i o{\u g}lu and Morris describe a possible extension of
   James' work for Weyl groups in general, where Young tableaux are 
   interpreted in terms of root systems. In this paper we show how to 
   construct submodules of Specht modules for Weyl groups.

\noindent
{\it 1991 Mathematics subject classification}: Primary 20C33 and 20F55, Secondary 
 22E45}

\renewcommand{\baselinestretch}{1.35}
\noindent
{\bf 1. Introduction}

The representation theory of the symmetric groups over fields of 
characteristic zero is well devoloped and documented. 
The original approach was due to G Frobenius and 
I Schur followed independently by A Young.
 In the 1930's, W Specht {\bf [8]} presented an
 alternative approach to give a full set of irreducible modules  
 now called Specht modules. In 1976, G D James, in a very important 
  paper {\bf [6]}, introduced an easy and 
 ingenious
 costruction of all the irreducible modules of the symmetric groups over 
 an arbitrary field which reduce to Specht modules in the case of
 characteristic zero. Al-Aamily,   
Morris and Peel {\bf [1]} showed how this 
construction could be extended to deal with the  Weyl groups of type 
$B_{n}$. 
In {\bf [7]} Morris described a possible extension of 
James' work for Weyl groups in general. 
An alternative and improved approach was described by the present author
 and Morris {\bf [5]}.  Later on the present author {\bf [4]} 
 develop the theory and
  show how a $K$-basis for Specht modules can be constructed in terms of
   standard tabloids and introduce an algorithm which results a basis for
    a Specht modules for Weyl groups in many cases. Recently, a possible
     extension of James' work for finite reflection groups has been given
      in {\bf [3]}. In this paper we further develop the theory and show 
    how to construct
 submodules of Specht modules for Weyl groups.

\noindent
{\bf 2. Notation and preliminary results}

In this section we establish the notation and state preliminary results on 
 Specht modules for Weyl groups. The basic definitions and background material required here may be found
  in {\bf [5]}. 

\noindent
{\bf 2.1} 
Let $V$ be $l$-dimensional Euclidean space over the real field  
{\bf R} equipped with a positive
 definite inner product ( , ). For $ \alpha \in V~,~\alpha \neq 0$ , let 
$\tau_{\alpha}$ be the $reflection$ in the hyperplane orthogonal to 
$\alpha$ , that is , $\tau_{\alpha}$ is the linear transformation on $V$ defined by 
\[ \tau_{\alpha}~(~v~)=v~-~2 ~\frac{(~\alpha~,~v~)}
{(~\alpha~,~\alpha~)}~ \alpha\]
\noindent
for all $ v \in V$ . Let  $\Phi$ be a root system
 in $V$ and $\pi$ be a simple
 system in $\Phi$ with corresponding positive system 
$\Phi^{+}$ and
 negative system $\Phi^{-}$ . Then , the $Weyl~group$ of $\Phi$
is the finite reflection group $ {\cal W}={\cal W}(\Phi)$ 
which is generated
by the $\tau_{\alpha},~\alpha \in \Phi$.

If $\Psi$  is a subsystem of $\Phi$ with simple
system $J \subset \Phi^{+}$ and Dynkin diagram $\Delta$ and let 
$\Psi=\displaystyle \bigcup_{i=1}^{r}\Psi_{i} $, where $\Psi_{i}$ are the
indecomposable components of $\Psi$. Let $J_{i}$ be a simple
system in $\Psi_{i}$ ($i=1,2,...,r$) and $J = \displaystyle
\bigcup_{i=1}^{r}J_{i}$. Let $\Psi^{\perp}$ be 
the largest subsystem
in $\Phi$ orthogonal to $\Psi$ and let $J^{\perp} 
\subset
\Phi ^{+}$ be the simple system of $\Psi^{\perp}$. 
If $\Psi^{'}$ is a subsystem of $\Phi$ which is
 contained in $\Phi
\setminus \Psi$, with simple system $J^{'} \subset 
\Phi^{+}$ and
Dynkin diagram $\Delta^{'}$,  $\Psi^{'}=\displaystyle
\bigcup_{i=1}^{s}\Psi_{i}^{'} $, where $\Psi_{i}^{'}$ 
are the
indecomposable components of $\Psi^{'}$ then let
 $J_{i}^{'}$ be a
simple system in $\Psi_{i}^{'}$ ($i=1,2,...,s$) and 
$J =
 \displaystyle
\bigcup_{i=1}^{s}J_{i}^{'}$.  Let $\Psi^{'^{\perp}}$ be the 
largest
subsystem in $\Phi$ orthogonal to $\Psi^{'}$ and let
 $J^{'^{\perp}}
\subset
\Phi ^{+}$ be the simple system of $\Psi^{'^{\perp}}$.
 If $\bar{J}$ stand for the ordered set $\{J_{1},J_{2},...,J_{r};
J_{1}^{'},J_{2}^{'},...,J_{s}^{'}\}$, where in addition 
the elements
in each $J_{i}$ and $J_{i}^{'}$ are also ordered, then 
let ${\cal T}_{\Delta}=\{ w\bar{J} \mid w\in {\cal W} \}$.
The pair $\bar{J} = \{J,J^{'}\}$ is called a $useful~ system$ 
in $\Phi$ if
${\cal W}(J)\cap {\cal W}(J^{'})=<e>$ and ${\cal W}(J^{\perp}) 
\cap
{\cal W}(J^{'^{\perp}})~=~~~<e>$. The
elements of ${\cal T}_{\Delta}$ are called $\Delta -tableaux$, 
the 
$J_{i}$ and $J_{i}^{'}$ are called the $rows$ and the $columns$ of 
$\{J,J^{'}\}$ respectively.
 Two $\Delta$-tableaux $\bar{J}$ and $\bar{K}$ are 
$row-equivalent$, written $\bar{J}~\sim~\bar{K}$, if
 there exists 
$w \in {\cal W}(J)$ such that $\bar{K}=w~\bar{J}$. 
 The equivalence class which contains the $\Delta$-tableau 
$\bar{J}$ is
denoted by  
$\{\bar{J}\}$ and is called a $\Delta-tabloid$. Let
 $\tau_{\Delta}$ be
set of all $\Delta$-tabloids. Then $\tau_{\Delta}=
\{ \{~\overline{dJ}~\} \mid d\in D_{\Psi} \}$, where 
$D_{\Psi}=\{~w \in {\cal W} \mid w~(j) \in \Phi^{+}~ for~
all~j \in J~\}$ is a distinguished set of coset representatives 
of ${\cal W}
(\Psi)$ in ${\cal W}$ . 
The group ${\cal W}$ acts on $\tau_{\Delta}$ as
$\sigma~\{\overline{wJ}\} =\{\overline{\sigma w J}\}$ 
for~all~$\sigma \in {\cal W}$.
Let $K$ be arbitrary field, let $M^{\Delta}$ be the
$K$-space whose basis elements are the $\Delta$-tabloids. 
 Extend the
action of ${\cal W}$ on $\tau_{\Delta}$ linearly on
 $M^{\Delta}$,
then $M^{\Delta}$ becomes a $K{\cal W}$-module. Let
\begin {eqnarray*}
\kappa_{J^{'}}=\sum_{\sigma~
\in~W(J^{'})}~s~(~\sigma~)~
\sigma~~~~~~and~~~e_{J,J^{'}}=\kappa_{J^{'}}~\{~\bar{J}~\}
\end{eqnarray*}
\noindent
where $s(\sigma)=(-1)^{l(\sigma)}$ is the sign function 
and $l(\sigma)$ is
the length of $\sigma$.  Then $e_{J,J^{'}}$
is called the generalized $\Delta- polytabloid $ associated with
$\bar{J}$. Let $S^{J,J^{'}}$ be the subspace of $M^{\Delta}$
generated by $ e_{wJ,wJ^{'}}~$ where $~w\in~{\cal W}$.  Then 
$S^{J,J^{'}}$
is called a $generalized~Specht~module$.
 A useful system $\{J,J^{'}\}$ in $\Phi$ is called a $good~system$ 
if $d~\Psi
\cap \Psi^{'} = \emptyset$ for $d \in D_{\Psi}$ then
$\{~\overline{dJ}~\}$ appears with non-zero coefficient in
$e_{J,J^{'}}$. If  $\{J,J^{'}\}$ is a good system ,
 then $S^{J,J^{'}}$ is  irreducible. {\bf [5]}

\noindent 
{\bf 2.2} 
Let $\Psi$ be a subsystem of $\Phi$ with simple system $J$.
If ${\cal N}(\Psi)~=~\{ w \in {\cal W} \mid w~\Psi=\Psi~\}$ and  
${\cal N}(J)~=~\{~w \in {\cal W} \mid w~(J)=J~\}$ , then 
${\cal N}(\Psi)~=~{\cal W}(\Psi)\propto {\cal N}(J)$, 
where $\propto$ denotes semidirect product {\bf [2]}.

\noindent
{\bf 3. Submodules of Specht Modules for Weyl Groups}

In this section we show how to construct submodules of Specht modules
 for Weyl groups. Detailed proofs are not always included as they are
  either given or modifications of those in the earlier paper {\bf [5]}.

Let $\bar{J}$ be the ordered set $\{J_{1},J_{2},...,
J_{r};J_{1}^{'},J_{2}^{'},...,J_{s}^{'}\}$  , 
satisfying (2.1) and let
\[{\cal T}_{\Delta}=\{ w\bar{J} \mid w\in {\cal W} \} \]

Now we can give our principal definition .

\noindent
{\bf Definition 3.1.} Let $\Psi$ and $\Psi^{'}$ be subsystems of $\Phi$ with 
simple systems $J$ and $J^{'}$ respectively such that 
$\Psi^{'} \subseteq \Phi \setminus \Psi$ and $J \subset \Phi^{+}$, 
$J^{'} \subset \Phi^{+}$. The pair $\{J,J^{'}\}$ is called a 
$useful ~sub-system $ in $\Phi$ if ${\cal N}(\Psi) \cap {\cal W}(\Psi^{'})=<e>$ and 
${\cal W}(\Psi^{\perp}) \cap {\cal W}(\Psi^{'^{\perp}})=
<e>$. 

\noindent
{\bf Remark 1.} If $\{J,J^{'}\}$ is a useful sub-system in $\Phi$ 
then $\{J,J^{'}\}$ is a useful system in $\Phi$.
If ${\cal N}(J)= <~e~>$, then the converse is true. 

\noindent
{\bf Definition 3.2.} Let $\{J,J^{'}\}$ be a useful sub-system  in $\Phi$. 
Then the elements of ${\cal T}_{\Delta}$ are called $\Delta-tableaux$.

\noindent
{\bf Definition 3.3.}  Two $ \Delta$-tableaux $\bar{J}$ and $\bar{K}$  
are $row-equivalent$ ,  written $\bar{J}~\sim~\bar{K}$ ,  if 
there exists  $w \in {\cal N}(\Psi)$ such that $\bar{K}~=~w~\bar{J}$.
The equivalence class which contains the $\Delta$-tableau $\bar{J}$ is 
$\{\bar{J}\}$ and is called a $\Delta^{*}$-tabloid.

Let $\tau_{\Delta^{*}}$ be the set of all $\Delta^{*}$-tabloids. 
It is clear that the number of distinct elements in 
$\tau_{\Delta^{*}}$ is  $[{\cal W}~:{\cal N}(\Psi)]$  and 
if $E_{\Psi}$ is the set of left coset representatives of 
${\cal N}(\Psi)$ in ${\cal W}$, then we have                               

\[\tau_{\Delta^{*}}=\{ \{~d\bar{J}~\} \mid d\in E_{\Psi} \} \]
\noindent 
We note that since ${\cal W}(\Psi)$ is a subgroup of ${\cal N}(\Psi)$, then 
$E_{\Psi} \subseteq D_{\Psi}$ and if $\bar{J}=\{~J~;~J^{'}~\}$, then 
$dJ \subset \Phi^{+}$ but $dJ^{'}$ need not be a subset of $\Phi^{+}$.

We now give an example to show how to construct a $\Delta^{*}$-tabloid. 
In this example and later examples we use the following notation. 
If $\pi = \{\alpha_{1},\alpha_{2},...,\alpha_{n}\}$ is a simple system in $\Phi$ and $\alpha \in \Phi$, then $\alpha= {\displaystyle \sum_{i=1}^{n}a_{i}\alpha_{i}}$, where $a_{i} \in {\bf Z}$. From now on $\alpha$ is 
denoted by $a_{1}a_{2}...a_{n}$ and $\tau_{\alpha_{1}},
\tau_{\alpha_{2}},...,\tau_{\alpha_{n}}$ are denoted by 
$\tau_{1},\tau_{2},...,\tau_{n}$ respectively. 

\noindent
{\bf Example 3.4.} Let $\Phi=\bf{A_{3}}$ with simple system  
$\pi=\{\alpha_{1}=\epsilon_{1}-\epsilon_{2},\alpha_{2}=\epsilon_{2}
-\epsilon_{3},
\alpha_{3}=\epsilon_{3}-\epsilon_{4}\}$.  
 Let $\Psi={\bf 2A_{1}}$ be
 the subsystem of $\Phi$ with simple system 
 $J = \{100,001\}$. Let $\Psi^{'}={\bf A_{1}}$ be the subsystem of $\Phi$  which is contained in 
$\Phi \setminus \Psi$, with simple system  $J^{'}=\{110\}$. 
Since ${\cal N}(\Psi) \cap {\cal W}(J^{'})=<e>$ and   
${\cal W}(\Psi^{\perp}) \cap {\cal W}(\Psi^{'^{\perp}})=<e>$, 
then $\{J,J^{'}\}$ is a useful sub-system in $\Phi$. 
Then $\tau_{\Delta^{*}}$ contains the $\Delta^{*}$-tabloids;

\begin {tabular}{llll}
$\{ \bar{J} \}$&=&$\{100,001;110\}~$&
\\
$\{ \tau_{2}\bar{J} \}$&=&$\{110,011;100 \}~$&
\\
$\{ \tau_{1}\tau_{2}\bar{J }\}$&=&$\{010,111;-100 \}~$&
\end{tabular}

\vspace{0.1cm}
The group ${\cal W}$ acts on $\tau_{\Delta^{*}}$ according to 
\[\sigma~\{\overline{wJ}\} =\{\overline{\sigma w J}\}~~~~~~for~all~\sigma \in {\cal W} .\]
\noindent
This action is easily seen to be well defined.

Now let $K$ be an arbitrary field and let $M^{\Psi}$ be the 
$K$-space whose basis elements are the $\Delta^{*}$-tabloids. 
Extend the action of ${\cal W}$ on $\tau_{\Delta}^{*}$ linearly on $M^{\Psi}$ , 
then $M^{\Psi}$ becomes a $K{\cal W}$-module. 
Then the~ $K{\cal W}$-module ~$M^{\Psi}$~  is~ a~ cyclic 
$K{\cal W}$-module~ generated~ by~ any~ one~ tabloid and 
~$dim_{K}M^{\Psi} = [{\cal W}:{\cal N}(\Psi)]$. 
Since $E_{\Psi} \subseteq D_{\Psi}$, then ~ $M^{\Psi}$~ is~ a~ $K{\cal W}$
-submodule~ of~ $M^{\Delta}$.

Now we proceed to consider the possibility of constructing a 
$K{\cal W}$-module which corresponds to the Specht module 
in the case of symmetric groups. In order to do this  we need to define 
a $\Delta^{*}$-polytabloid .

\noindent
{\bf Definition 3.5.} Let $\{J,J^{'}\}$ be a useful sub-system in $\Phi$. 
Let
\begin {eqnarray*}
\kappa_{J^{'}}~=~\sum_{\sigma~\in~{\cal W}(\Psi^{'})}~s~(~\sigma~)~\sigma~~~~~~and~~~e_{J,J^{'}}~=~\kappa_{J^{'}}~\{~\bar{J}~\}
\end{eqnarray*}
\noindent
where $s$ is the sign function. 
Then  $e_{J,J^{'}}$ is called the generalized $\Delta^{*}-polytabloid$.

If $w \in {\cal W}$, then $w~e_{J,J^{'}}=\kappa_{wJ^{'}}~\{~\overline{wJ}~\}
~=~e_{wJ,wJ^{'}}$. Hence, if $S^{\Psi,\Psi^{'}}$ is the subspace of $M^{\Psi}$ 
 generated by $e_{wJ,wJ^{'}}$ where $w \in {\cal W}$, then 
 $S^{\Psi,\Psi^{'}}$ is a 
$K{\cal W}$-submodule of $M^{\Psi}$ and $S^{\Delta,\Delta^{'}}$, 
which is called a $submodule~of$ $generalized$ $Specht$ $module$. The~ $K{\cal W}$-module~ 
$S^{\Psi,\Psi^{'}}$~ is~ a 
~cyclic~ submodule~ generated~ by~ any $\Delta^{*}-polytabloid$.

The submodule of Specht module is spanned by the $e_{wJ,wJ^{'}}$ 
for all $ w \in {\cal W} $. But it is easy to show that 
$S^{\Psi,\Psi^{'}}$~ is~ generated~ by $e_{dJ,dJ^{'}}$, 
where~ $d \in D_{\Psi^{'}}$. Now we can give the following lemma.

\noindent
{\bf Lemma 3.6.} $Let ~\{J,J^{'}\}~ be~ a~ useful~ sub-system ~in~\Phi 
~and~ 
let ~d \in E_{\Psi}.~ If~ \{\overline{dJ} \}~ appears$ 
\newline
$in~ e_{J,J^{'}}~ 
~then ~it~ appears~ only~ once.$

\noindent
{\bf Proof.} If $\sigma~,~\sigma^{'} \in {\cal W}(\Psi^{'})$ and 
suppose that $\sigma = dw~,~
\sigma^{'}=dw^{'}$ where $w,w^{'} \in {\cal N}(\Psi) $. 
Then $d=\sigma w^{-1}=
\sigma^{'}w^{'^{-1}}$ and $\sigma^{'^{-1}}\sigma=w^{'^{-1}}w \in 
{\cal N}(\Psi) 
\cap {\cal W}(\Psi^{'})=~<~e~>$. 
Hence we have $w=w^{'}$ and $\sigma=\sigma^{'}$. 
Then $\{\overline{dJ} \}$ appears in $e_{J,J^{'}}$ only once.$~~~~\Box$

\noindent
{\bf Corollary 3.7.} $If~\{J,J^{'}\}~ is~ a~ 
useful~ sub-system~ in~\Phi,~ then ~e_{J,J^{'}} \neq 0$.

\noindent
{\bf Proof.} By Lemma 3.6 if $\{\overline{\sigma J} \}$ appears in 
$e_{J,J^{'}}$ 
then  all the $\{\overline{\sigma J }\}$ are different, 
where $\sigma \in {\cal W}(\Psi^{'})$.  
But $\{~\{\overline{\sigma J }\} \mid \sigma \in {\cal W}(\Psi^{'})~\}$ 
is a linearly 
independent subset of $\{~\{\overline{d~J }\} \mid d \in E_{\Psi}~\}$.   
If $e_{J,J^{'}} =0$ then $s(\sigma)=(-1)^{l(\sigma)}=0$ for all $\sigma 
\in {\cal W}(\Psi^{'})$ . This 
is a contradiction and so  $e_{J,J^{'}} \neq 0.$ $~~~\Box$ 

\noindent
{\bf Remark 2.} In {\bf [5]} the present author and Morris defined a 
useful system as a pair $\{J,J^{'}\}$ such that 
 ${\cal W}(\Psi) \cap {\cal W}(\Psi^{'})=<e>$ and 
${\cal W}(\Psi^{\perp}) \cap {\cal W}(\Psi^{'^{\perp}})=
<e>$. However Lemma 3.6 and the following lemma show that in this 
new setting the condition
${\cal N}(\Psi) \cap {\cal W}(\Psi^{'})=<e>$ in our definition of 
a useful sub-system is necessary. 

\noindent
{\bf Lemma 3.8.} $If~ there~ exists ~w \in {\cal N}(\Psi) \cap 
{\cal W}(\Psi^{'})~ such~ that~w~ has~ order~ 2, 
~and~ s(w)~ = -1$  
\newline
$then ~e_{J,J^{'}}~=~0$.

\noindent
{\bf Proof.}  See Lemma 3.14 {\bf [5]}.$~~~~\Box$

\noindent
{\bf Example 3.9.} Let $\Phi={\bf G_{2}}$ and 
$\pi=\{\alpha_{1}=\epsilon_{1}-\epsilon_{2},\alpha_{2}
=-2\epsilon_{1}+\epsilon_{2}+\epsilon_{3}\}$. Let $\Psi=
{\bf A_{1}}$ be the subsystem of $\Phi$ with simple system 
$J=\{\alpha_{1}\}$ and 
let $\Psi^{'}={\bf A_{2}}$ be the subsystem of $\Phi$ with simple system 
$J^{'}=\{\alpha_{2},3\alpha_{1}+\alpha_{2}\}$. Thus we have 
${\cal W}(\Psi)=\{e,\tau_{1}\}$  and
\newline 
${\cal W}(\Psi^{'})=\{e,\tau_{2},\tau_{1}\tau_{2}\tau_{1},
\tau_{2}\tau_{1}\tau_{2}\tau_{1},\tau_{1}\tau_{2}\tau_{1}\tau_{2},
\tau_{2}\tau_{1}\tau_{2}\tau_{1}\tau_{2}\},~{\cal N}(\Psi)=~\{e,\tau_{1},
\tau_{2}\tau_{1}\tau_{2}\tau_{1}\tau_{2},
\tau_{2}\tau_{1}\tau_{2}\tau_{1}\tau_{2}\tau_{1}\}$. 
Then ${\cal W}(\Psi) \cap {\cal W}(\Psi^{'})=<e>$, but 
$w=~\tau_{2}\tau_{1}\tau_{2}\tau_{1}\tau_{2} \in 
{\cal N}(\Psi) \cap {\cal W}(\Psi^{'}) $ and $e_{J,J^{'}}=0$.

If  $\{J,J_{1}^{'}\}$ and $\{J,J_{2}^{'}\}$ are~ 
sub-useful~ systems~ in~$\Phi$ and $\Psi_{1}^{'} \subseteq \Psi_{2}^{'}$, 
then $S^{\Psi,\Psi_{2}^{'}}$~is~ a~ $K{\cal W}$-submodule~ of~ 
$S^{\Psi,\Psi_{1}^{'}}$,~where $J_{1}^{'}$~ and $J_{2}^{'}$ 
are~ simple~ systems~ for $\Psi_{1}^{'}$ and $\Psi_{2}^{'}$ respectively.

Now we  consider under what conditions  $S^{\Psi,\Psi^{'}}$ is irreducible .

\noindent
{\bf Lemma 3.10.} $Let~ \{J,J^{'}\}~ be~ a~ useful~ sub-system~ in~\Phi~
and~  let~ d \in E_{\Psi}.~ Then~ the$  
\newline
$following~ conditions~ are ~equivalent:$
\newline
$(i)~ \{~\overline{dJ}~\}~ appears~ with~ non-zero~ coefficient~ in~
e_{J,J^{'}}$
\newline
$(ii)~ There~ exists~\sigma \in {\cal W}(\Psi^{'})~such~ that 
~\sigma \{~\bar{J}~\}~=~\{\overline{dJ}~\}$
\newline
$(iii)~ There~ exists~\rho \in {\cal N}(\Psi)~  and ~\sigma \in 
{\cal W}(\Psi^{'})~ such ~that~d~=~\sigma~\rho.$

\noindent
{\bf Proof.} This follows from Lemma 3.18 {\bf [5]} $~~~~~~~~\Box $

\noindent
{\bf Lemma 3.11.} $Let ~\{J,J^{'}\}~ be~ a~ useful~ sub-system~ in~
\Phi~ and ~ let 
~d \in E_{\Psi}.~ If ~\{~\overline{dJ}~\}~ appears$ 
 $in~e_{J,J^{'}}~ then~d~\Psi \cap \Psi^{'}~=~\emptyset $.

\noindent
{\bf Proof.}
 Let $\alpha \in d~\Psi$. If $\{~\overline{dJ}~\}$ appears 
in $e_{J,J^{'}}$ then by Lemma 3.10 $d~=~\sigma~\rho$, where $\sigma \in 
{\cal W}(\Psi^{'})$ and $\rho \in {\cal N}(\Psi)$. Then 
$ \alpha \in \sigma \rho \Psi$. Since $\rho \in {\cal N}(\Psi)$, 
then $ \alpha \in \sigma \Psi$ and $\sigma^{-1}\alpha \in \Psi$. But 
$\Psi \cap \Psi^{'} = \emptyset$, then $\sigma^{-1}\alpha \not \in \Psi^{'}$. 
Since $\sigma \in {\cal W}(J^{'})$, $\sigma \Psi^{'}~=~\Psi^{'}$ 
then $\alpha \not \in \Psi^{'}$. $~~~~\Box$

\noindent
{\bf Lemma 3.12.} $Let~\{J,J^{'}\}~  be~  a~ useful~ sub-system~ in~\Phi~
and 
~let~ d \in E_{\Psi}.~ Let~d~\Psi \cap \Psi^{'} \not = \emptyset$. 
 $Then ~\kappa_{J^{'}}\{~\overline{dJ}~\} = 0 $.

\noindent
{\bf Proof.} If $d~\Psi \cap \Psi^{'} \not = \emptyset$ then let 
$\alpha \in d~\Psi \cap \Psi^{'} $ and so $\tau_{\alpha} \in
 {\cal N}(d\Psi) \cap {\cal W}(\Psi^{'})$. Thus 
\[(~e-\tau_{\alpha}~)~\{\overline{dJ}\} =
 \{\overline{dJ}\} - \{\overline{dJ}\} = 0 \]
\noindent
Since $\{~e~,~\tau_{\alpha}~\}$ is a 
subgroup of ${\cal W}(\Psi^{'})$ then we can select signed coset 
representatives $\sigma_{1}$, $~\sigma_{2}$, ~$\sigma_{3}$, ~...~, ~$\sigma_{s}$ for $\{~e~,~\tau_{\alpha}~\}$ in ${\cal W}(J^{'})$ such that

\[\kappa_{J^{'}}\{~\overline{dJ}~\} =(~\sum_{i~=~1}^{s}~\sigma_{i})~(~e~-~\tau_{\alpha}~)~\{\overline{d~J}\}~=~0 ~~. ~~~~~~~~\Box \]

The converse of Lemma 3.11 is not true in general, which leads to 
the following definition .

\noindent
{\bf Definition 3.13.}  A useful sub-system $\{J,J^{'}\}$ in $\Phi$ is 
called a $good~sub-system$  if 
$d\Psi \cap \Psi^{'} = \emptyset$ for $d \in E_{\Psi}$ 
then $\{~\overline{dJ}~\}$ 
appears with non-zero coefficient in $e_{J,J^{'}}$.

\noindent
{\bf Lemma 3.14.} $Let~\{J,J^{'}\}~ be~ a ~good ~sub-system~ in~\Phi~ and ~ 
let ~d \in E_{\Psi}$.
\newline
$(i)~ If ~\{\overline{dJ}\}~ does~ not ~appear~ in~e_{J,J^{'}}~ 
then ~\kappa_{J^{'}}\{~\overline{dJ}~\} = 0 $.
\newline
$(ii)~ If ~\{\overline{dJ}\}~  appears~ in ~e_{J,J^{'}}~  
then ~there ~exists ~\sigma \in {\cal W}(\Psi^{'})~ such~ that$
\[\kappa_{J^{'}}\{~\overline{dJ}~\}~=~s~(~\sigma~)~e_{J,J^{'}} \]

\noindent
{\bf Proof}  See Lemma 3.22 {\bf [5]}. $~~~~~~~\Box$

By the previous lemma if $m \in M^{\Psi}$, then $\kappa_{J^{'}}m$ 
is~ a multiple of $e_{J,J^{'}}$.

We now define a bilinear form $ <~,~>$ on 
$M^{\Psi}$ by setting

\[
<~\{\bar{J_{1}}\}~,~\{\bar{J_{2}}\}~>~=
\left\{
\begin{array}{ll}
1& \mbox{if~$\{\bar{J_{1}}\}~=~\{\bar{J_{2}}\}$}\\
0& \mbox{otherwise} 
\end{array}
\right.
\]

\noindent
This is a symmetric, non-singular, ${\cal W}$-invariant bilinear form
 on $M^{\Psi}$ .

Now we can give James' submodule theorem in this new setting.

\noindent
{\bf Theorem 3.15.} $Let~\{J,J^{'}\}~ be~ a~ good~sub-system~ in ~\Phi.~
Let~ U ~be~ submodule~ of~ M^{\Psi}$.
\newline 
$Then~ either ~S^{\Psi , \Psi^{'}} \subseteq U~or~U 
\subseteq S^{\Psi,
 \Psi^{'^\perp}},~ where ~S^{\Psi,\Psi^{'^\perp}}~is~ the~complement~ of~ 
 S^{\Psi,\Psi^{'}}~ in~M^{\Psi}.$

\noindent
{\bf Proof.} See Theorem 3.24 {\bf [5]}.$~~~~~~~\Box$

We can now prove our principal result.

\noindent
{\bf Theorem 3.16.} $Let~\{J,J^{'}\}~be~ a~ good~ sub-system~ in~\Phi.~ 
The~K{\cal W}-module$
\[
D^{\Psi,\Psi^{'}}~=~S^{\Psi,\Psi^{'}}~
/~S^{\Psi,\Psi^{'}}\cap S^{\Psi,\Psi^{'^\perp}}\] 
$is~ zero~ or~ irreducible.$

\noindent
{\bf Proof.} 
 If $U$ is a submodule of $S^{\Psi,\Psi^{'}}$ 
then $U$ is a 
submodule of $M^{\Psi}$ and by Theorem 3.15 
either $S^{\Psi,\Psi^{'}} 
\subseteq U$ in which case $U~=~
S^{\Psi,\Psi^{'}} $ or $U \subseteq 
S^{\Psi,\Psi^{'^\perp}}$ and $U
 \subseteq ~S^{\Psi,\Psi^{'}}~\cap ~
S^{\Psi,\Psi^{'^\perp}}$ , which 
completes the proof. $~~~~~\Box$

\noindent
{\bf Remark 3.} 
Since $E_{\Psi} \subseteq D_{\Psi}$, if $\{J,J^{'}\}$ is a useful
 sub-system and $S^{\Delta,\Delta^{'}}$ is 
 irreducible, then $S^{\Psi,\Psi^{'}}$ is irreducible.

In {\bf [5]}, we see that not all the irreducible modules for 
${\cal W}({\bf D_{4}})$  are obtained by considering subsystems of root
 systems. In order to obtain additional irreducible character, we define
extended subsystems of root systems. It is shown that 
  the irreducible character 
  $\chi_{12}$ of degree 6 is obtained by considering 
  ${\cal W}({\bf G_{2}})$ as Steinberg subgroup of ${\cal W}({\bf D_{4}})$.

  In the following example we show how the irreducible character 
  $\chi_{12}$ of degree 6 can be obtained by considering a good sub-system 
  and how a  good sub-system 
may be constructed in the Weyl group of 
type ${\bf D_{4}}$.

\noindent
{\bf Example 3.17.} Let $\Phi={\bf D_{4}}$ 
 with simple system 
 \[\pi =\{\alpha_{1}=\epsilon_{1}-\epsilon_{2},
\alpha_{2}=\epsilon_{2}-\epsilon_{3},\alpha_{3}=\epsilon_{3}-\epsilon_{4},
\alpha_{4}=\epsilon_{3}+\epsilon_{4}\}\] 
\noindent
The representatives
 of conjugate classes and the character table  of 
 ${\cal W}({\bf D_{4}})$ and the 
 non-conjugate subsystems of ${\bf D_{4}}$ have been given in {\bf [5]}.

Let $\Psi={\bf A_{3}^{'}}$ be the subsystem 
of ${\bf D_{4}}$ with $J=\{1000,0100,0001\}$. 
Let $\Psi^{'}={\bf A_{1}}$ 
 be the subsystem of $\Phi$ which is 
contained in $\Phi \setminus \Psi$, with simple system 
$J^{'}=\{1110\}$. Then ${\cal N}(\Psi)= {\cal W}({\bf B_{3}})$ and 
since ${\cal N}(\Psi) \cap 
{\cal W}(\Psi^{'})=<e>$ and   
${\cal W}(\Psi^{\perp}) \cap {\cal W}(\Psi^{'^{\perp}})=<e>$, 
then $\{J,J^{'}\}$ is a useful sub-system in $\Phi$. 
Then $\tau_{\Delta^{*}}$ contains the $\Delta^{*}$-tabloids;

\begin {tabular}{lll}
$\{ \bar{J} \}$&=&$\{1000,0100,0001;1110\}~$\\
$\{ \tau_{3}\bar{J} \}$&=&$\{1000,0110,0001;1100 \}~$\\
$\{ \tau_{2}\tau_{3}\bar{J }\}$&=&$\{1100,0010,0101;1000\}~$\\
$\{ \tau_{1}\tau_{2}\tau_{3}\bar{J}\}$&=
&$\{0100,0010,1101;-1000\}~$
\end{tabular}

\vspace{0.1cm}
\noindent
For $d=e,\tau_{1}\tau_{2}\tau_{3}$ we have  
$d\Psi\cap\Psi^{'}=\emptyset$. 
Since
\[
e_{J,J^{'}}=
\{\bar{J} \}-\{ \tau_{1}\tau_{2}\tau_{3}\bar{J} \} 
\]
\noindent
then $\{J,J^{'}\}$ is a good sub-system in $\Phi$ .

Now  let $K$ be a field with Char$K$=0. Let $M^{\Psi}$
 be $K$-space whose basis elements are the $\Delta^{*}$-tabloids. 
 Let $S^{\Psi,\Psi^{'}}$ be the corresponding 
 $K{\cal W}$-submodule of 
 $M^{\Psi}$, then  by definition of 
 the submodule of Specht module we have 
\[S^{\Psi,\Psi^{'}}=
Sp~\{~e_{J,J^{'}}~,~e_{\tau_{3}J,\tau_{3}J^{'}}~,~
e_{\tau_{2}\tau_{3}J,
\tau_{2}\tau_{3}J^{'}}~\}\]

\noindent
where

$\begin{array}{lll}
e_{J,J^{'}}&=&\{\bar{J} \}-
\{ \tau_{1}\tau_{2}\tau_{3}\bar{J} \}
\\
e_{\tau_{3}J,\tau_{3}J^{'}}&
=&\{\tau_{3}\bar{J} \}-\{ \tau_{1}\tau_{2}
\tau_{3}\bar{J} \}
\\
e_{\tau_{2}\tau_{3}J,\tau_{2}\tau_{3}J^{'}}&=
&\{\tau_{2}\tau_{3}
\bar{J} \}-
\{ \tau_{1}\tau_{2}
\tau_{3}\bar{J}\}
\end{array}$

Let $T$ be the matrix representation of ${\cal W}$ afforded by 
$S^{\Psi,\Psi^{'}}$ with character $\psi$  and 
let $\tau_{1}\tau_{3}\tau_{2}$ be the representative of the
conjugate class $C_{4}$ as in {\bf [5]}. Then

{\setlength{\arraycolsep}{0.1cm}
$\begin{array}{lllll}
\tau_{1}\tau_{3}\tau_{2}(e_{J,J^{'}})&=
&\{\tau_{3}\bar{J} \}-\{ \tau_{2}
\tau_{3}\bar{J} \}
&=&
e_{\tau_{3}J,\tau_{3}J^{'}}~-
~e_{\tau_{2}\tau_{3}J,\tau_{2}\tau_{3}J^{'}}\\

\tau_{1}\tau_{3}\tau_{2}(e_{\tau_{3}J,\tau_{3}J^{'}})&
=&\{ \tau_{1}\tau_{2}
\tau_{3}\bar{J}\}-\{ \tau_{2}\tau_{3}\bar{J} \}
&=&
-~e_{\tau_{2}\tau_{3}J,\tau_{2}\tau_{3}J^{'}}\\

\tau_{1}\tau_{3}\tau_{2}(e_{\tau_{2}\tau_{3}J,\tau_{2}\tau_{3}J^{'}})
&=&
\{\bar{J} \}-\{ \tau_{2}
\tau_{3}\bar{J}\}
&=&
e_{J,J^{'}}~-~ e_{\tau_{3}J,\tau_{3}J^{'}}
\end{array}$}

Thus we have

$T~(~\tau_{1}\tau_{3}\tau_{2}~)~=
~\left( \begin{array}{ccc}
0&1&-1\\
0&0&-1\\
1&0&-1
\end{array}
\right) $
     and $\psi(\tau_{1}\tau_{3}\tau_{2})=-1$.

By a similar calculation to the above it can be showed that $\psi=
\chi_{9}$. Now let $\Psi={\bf A_{2}}$ be the subsystem 
of ${\bf D_{4}}$ with $J=\{1000,0100\}$. 
Let $\Psi^{'}={\bf A_{2}}$ 
 be the subsystem of $\Phi$ which is 
contained in $\Phi \setminus \Psi$, with simple system 
$J^{'}=\{0001,0110\}$. Then ${\cal N}(\Psi)= {\cal W}({\bf G_{2}})$,  
the corresponding character of ${\cal W}$ afforded by 
$S^{\Psi,\Psi^{'}}$ is $\chi_{12}$. 

\begin{center}        
REFERENCES
\end{center}

{\bf 1.} {\sc E. Al-Aamily , A . O .  Morris} and {\sc M . H . Peel},
The~ representations~ of~ the~ Weyl groups~ of ~type ~$B_{n}$, 
$Journal~ of~ Algebra$ {\bf 68} ( 1981), 298--305.

{\bf 2.} {\sc R . W . Carter},
Conjugacy~ classes~ in ~the Weyl~ group,
$Comp.~  Math.$, {\bf 25} (1972), 1--59.  

{\bf 3.} {\sc S. Hal\i c\i o{\u g}lu}, 
Specht Modules for Finite Reflection  Groups, 
appears in $Glasgow$ $Mathematical$ $Journal$.

{\bf 4.} {\sc S. Hal\i c\i o{\u g}lu}, 
 A Basis of Specht~ Modules~ for~ Weyl~ Groups, 
submitted~ for~ publication.

{\bf 5.} {\sc S. Hal\i c\i o{\u g}lu} and {\sc A. O. Morris}, 
Specht~ Modules~ for~ Weyl~ Groups, \\
$Contributions$ $to~ Algebra$ $and~ Geometry$, {\bf 34} (1993), 257--276.

{\bf 6.} {\sc G. D. James},
The~ irreducible~ representations~ of~ the~ symmetric~ group, \\
$Bull.~ Lond.~ Math.~ Soc$, {\bf 8} (1976), 229--232.

{\bf 7.} {\sc A. O. Morris},
Representations~ of~ Weyl~ groups~ over~ an~ arbitrary~ field, 
$Ast\grave{e}risque$ {\bf 87-88} (1981), 267--287.

{\bf 8.} {\sc W . Specht}, Die irreduziblen darstellungen der 
symmetrichen gruppe,  $Math~Z.$ {\bf 39} (1935), 679--711.

\noindent
{\sc Department of Mathematics\\ 
Ankara University\\ 
06100~  Tando{\u g}an\\
Ankara Turkey}\\

\end{document}